\documentclass[reqno]{amsart}

\calclayout

\usepackage[unicode]{hyperref}
\usepackage[lite]{amsrefs}
\usepackage{graphicx}

\def\paragraph{\subsection}
\makeatletter
\def\@sect#1#2#3#4#5#6[#7]#8{%
  \edef\@toclevel{\ifnum#2=\@m 0\else\number#2\fi}%
  \ifnum #2>\c@secnumdepth \let\@secnumber\@empty
  \else \@xp\let\@xp\@secnumber\csname the#1\endcsname\fi
  \@tempskipa #5\relax
  \ifnum #2>\c@secnumdepth
    \let\@svsec\@empty
  \else
    \refstepcounter{#1}%
    \edef\@secnumpunct{%
      \ifdim\@tempskipa>\z@ 
        \@ifnotempty{#8}{.\@nx\enspace}%
      \else
        \@ifempty{#8}{.}{.\@nx\enspace}%
      \fi
    }%
      \ifnum #2=\tw@ \def\@secnumfont{\bfseries}\fi{}%
    \protected@edef\@svsec{%
      \ifnum#2<\@m
        \@ifundefined{#1name}{}{%
          \ignorespaces\csname #1name\endcsname\space
        }%
      \fi
      \@seccntformat{#1}%
    }%
  \fi
  \ifdim \@tempskipa>\z@ 
    \begingroup #6\relax
    \@hangfrom{\hskip #3\relax\@svsec}{\interlinepenalty\@M #8\par}%
    \endgroup
    \ifnum#2>\@m \else \@tocwrite{#1}{#8}\fi
  \else
  \def\@svsechd{#6\hskip #3\@svsec
    \@ifnotempty{#8}{\ignorespaces#8\unskip
       \@addpunct.}%
    \ifnum#2>\@m \else \@tocwrite{#1}{#8}\fi
  }%
  \fi
  \global\@nobreaktrue
  \@xsect{#5}}
\makeatother

\makeatletter
\def\pxspace{\@ifnextchar.{\@}{.\@\xspace}}
\makeatother

\newcommand{\interject}[1]{\noalign{\begin{quote}#1\end{quote}}}

\newcommand{\disconnect}{\leavevmode\par}

\newcommand{\mbbN}{\mathbb{N}}
\newcommand{\mbbY}{\mathbb{Y}}

\newcommand{\Ideal}[1]{\texorpdfstring{$\{#1\}$}{\{#1\}}}

\begin{document}

\newcommand{\ipeFigEmpty}{1}
\newcommand{\ipeFigEmptyplusAB}{2}
\newcommand{\ipeFigAplusB}{3}
\newcommand{\ipeFigAplusBCD}{4}
\newcommand{\ipeFigBplusA}{5}
\newcommand{\ipeFigBplusAEF}{6}
\newcommand{\ipeFigABplusCDEF}{7}
\newcommand{\ipeFigACplusBD}{8}
\newcommand{\ipeFigACplusBDG}{9}
\newcommand{\ipeFigADplusBC}{10}
\newcommand{\ipeFigADplusBCH}{11}
\newcommand{\ipeFigACDplusBGH}{12}
\newcommand{\ipeFigACDplusBGHI}{13}
\newcommand{\ipeFigABCplusDEFG}{14}
\newcommand{\ipeFigBEplusAF}{15}
\newcommand{\ipeFigBEplusAFJ}{16}
\newcommand{\ipeFigBFplusAE}{17}
\newcommand{\ipeFigBFplusAEK}{18}
\newcommand{\ipeFigBEFplusAJK}{19}
\newcommand{\ipeFigBEFplusAJKL}{20}
\newcommand{\ipeFigP}{21}
\newcommand{\ipeFigPA}{22}
\newcommand{\ipeFigPAiB}{23}
\newcommand{\ipeFigPAiAjBC}{24}

\title[Differential distributed lattices]%
      {Constructing/analyzing differential distributed lattices}
\author{Dale R. Worley}
\email{worley@alum.mit.edu}
\date{Apr 12, 2026} 

\begin{abstract}
We restate a process presented by Stanley as a technique to prove that
there exists exactly one $d$-differential distributive lattice for any
positive integer $d$.  This process can be trivially extended to apply
to distributive finitary lattices that have a variety of differential
poset structures.  It can be viewed as an algorithm for constructing
such lattices.  Alternatively, it can be viewed as an algorithm for
analyzing and characterizing such lattices.

We show that the process can be used to prove properties of all
weighted-differential lattices with positive weights.  We present this
with the hope that this approach can be used as the basis for a
complete characterization of distributive lattices with a
weighted-differential structure with positive weights.
\end{abstract}

\maketitle

\textit{The 2\textsuperscript{nd} version corrects typos and improves
the figures.}
\textit{The 3\textsuperscript{rd} version adds a reference to the open
problem notebook \cite{Wor2026a}.}


\section{Introduction}

This article is a restatement of a process described by Stanley in
\cite{Stan1975a}*{sec.~3} and in the
solution of exercise~3.51(a) of \cite{Stan2012a},
which is presented as a technique to prove that there exists exactly
one $d$-differential distributive lattice for any positive integer $d$.
This process can be trivially extended to apply to distributive
finitary lattices that have a variety of differential poset structures.
It can be viewed as an algorithm for constructing such lattices,
which may be either deterministic or nondeterministic, depending on the
variety of differential structure under consideration.
Alternatively, it can be viewed as an algorithm for analyzing and
characterizing such lattices.

We present this with the hope that it can be used as the basis for a
complete characterization of distributive lattices with a
weighted-differential structure with positive weights.
Specifically, we show that the process can be used to prove
properties of all weighted-differential lattices with
positive weights.
Since these properties resemble the assumptions used in
the characterization \cite{Wor2026b}, it may be that the conclusions
of \cite{Wor2026b} apply with weaker hypotheses, ideally applying to
all distributive lattices with positive weights and a positive, constant
distributive degree.
This goal is problem 50 in \cite{Wor2026a}, the notebook of open
problems in the combinatorics of tableaux.

\section{The process}

We assume that the lattice is a finitary distributive lattice.
Thus it is isomorphic
to the lattice of finite ideals of the poset of join-irreducible
elements of the lattice.  The elements of the poset are called
\emph{points}.  We identify the elements of the lattice with the
corresponding ideal of points, which is the set of points $\leq$ the
given element.

The maximal points of an ideal are its \emph{deletion points}, the
points that can be deleted from it to produce a smaller ideal.

The minimal points that are \emph{not} in an ideal are its
\emph{insertion points}, the points that can be added to it to
produce a larger ideal.

The deletion points of an ideal correspond to elements of the lattice
covered by the corresponding lattice point.  In our figures,
deletion points are shown in blue.
The insertion points
correspond to the elements of the lattice that cover the corresponding
lattice point.
In our figures, insertion points are shown in orange.

Each point is
assigned a \emph{weight}, which is a positive integer.  Every lattice
element has a
\emph{differential degree}, which is a positive integer.  The lattice
is \emph{weighted differential} iff for every ideal,
$$ \sum \textup{weight of insertion point} =
\sum \textup{weight of deletion point} +
\textup{element's differential degree} $$

The process is defined in terms of ideals of points.
The process iterates through the the ideals of points in increasing
order, starting with the empty ideal, $\{\}$.
This corresponds to examining the elements of the lattice in
increasing order.  Each iteration consists of these steps:

\begin{enumerate}

\item Consider the ideal.

\item Determine its set of deletion points and the sum of their weights.

\item Calculate the sum of the weights of the ideal's insertion
  points:
\begin{equation*}
\textup{total weight of insertion points} =
  \sum \textup{weight of deletion point} +
  \textup{element's differential degree}
\end{equation*}

\item Determine the set of its \emph{old} insertion points, those that
  have been
  constructed/analyzed during iterations for smaller ideals.  These
  insertion points are necessarily not in
  the ideal.  They were seen in previous iterations, so none of them
  are above \emph{all} elements of the ideal.

\item Calculate the sum of the weights of the
  \emph{new} insertion points that have to be constructed/analyzed
  during this iteration:
$$ \textup{total weight of new insertion points} =
  \textup{total weight of insertion points} -
  \sum \textup{weight of old insertion point} $$

\item The total weight of new insertion points must be non-negative.
  If we are using the process to construct a lattice and the total
  weight is negative, the construction
  has failed.  If we are analyzing a given lattice, its existence proves
  this total weight is non-negative.

\item Add zero of more new insertion points whose weights total to
  \textit{total weight of new insertion points}.  If the weights are
  fixed at 1, there is only once choice, to add as many new insertion points
  as the total.  But if the weights can vary,
  this step may be nondeterministic.  Since the new insertion points must
  not have been insertion points seen in any previous iteration, they must each
  be greater than exactly all the points in the ideal.  Since we are
  drawing the Hasse diagrams of the ideals, each new insertion point
  covers exactly the deletion points of the ideal.

\item The existence of the new insertion points implies the existence
  of larger ideals whose existence has not been determined before (and
  which will be examined in later iterations).
  Specifically, for every known ideal containing this ideal and for
  every non-empty subset of the new insertion points, there is an
  ideal which is the union of that ideal and that subset.

\end{enumerate}

\section{Example}

A simple example of the process is the demonstration that there
exists exactly one $2$-differential distributive lattice.
That is, we will construct a weighted-differential distributive
lattice where the weight of all points is 1 and the differential
degree of all lattice elements is 2.
This construction succeeds and is deterministic, showing that the only
such lattice is $\mbbY^{\times 2}$, the cartesian square of Young's
lattice $\mbbY$ of partitions.
The poset of points of Young's lattice is the quadrant of
integral points in the plane, $\mbbN \times \mbbN$.  Thus the poset we
construct will be the disjoint union of two copies of the quadrant.

We will draw ideals of points in the usual ``upward is greater''
style of Hasse diagrams, so we show the quadrant in ``Russian
format'', with $(0,0)$ at the bottom-center.

\paragraph{The ideal \Ideal{}} \disconnect

\newcommand{\Figure}[1]{
  \vspace{\baselineskip}
  \begin{center}
  \includegraphics[page=#1]{dist-diff-figs.pdf}
  \end{center}
  \vspace{\baselineskip}}

\Figure{\ipeFigEmpty}

The deletion points are $\{\}$.
The number of insertion points is 2.
The existing insertion points are $\{\}$.
The new insertion points are $\{ A, B \}$.

\Figure{\ipeFigEmptyplusAB}

Thus there are future ideals $\{A\}$, $\{B\}$, and $\{A,B\}$.

\paragraph{The ideal \Ideal{A}} \disconnect

\Figure{\ipeFigAplusB}

The deletion points are $\{A\}$.
The number of insertion points is 3.
The existing insertion points are $\{B\}$.
The new insertion points are $\{C,D\}$.

\Figure{\ipeFigAplusBCD}

Thus there are future ideals $\{A,C\}$, $\{A,D\}$, $\{A,C,D\}$,
$\{A,B,C\}$, $\{A,B,D\}$, and $\{A,B,C,D\}$.

\paragraph{The ideal \Ideal{B}} \disconnect

\Figure{\ipeFigBplusA}

The deletion points are $\{B\}$.
The number of insertion points is 3.
The existing insertion points are $\{A\}$.
The new insertion points are $\{E,F\}$.

\Figure{\ipeFigBplusAEF}

Thus there are future ideals $\{B,E\}$, $\{B,F\}$, $\{B,E,F\}$,
$\{A,B,E\}$, $\{A,B,F\}$, $\{A,B,E,F\}$,
$\{A,B,C,E\}$, $\{A,B,C,F\}$, $\{A,B,C,E,F\}$,
$\{A,B,D,E\}$, $\{A,B,D,F\}$, $\{A,B,D,E,F\}$,
$\{A,B,C,D,E\}$, $\{A,B,C,D,F\}$, and $\{A,B,C,D,E,F\}$.

\paragraph{The ideal \Ideal{A,B}} \disconnect

\Figure{\ipeFigABplusCDEF}

The deletion points are $\{A,B\}$.
The number of insertion points is 4.
The existing insertion points are $\{C,D,E,F\}$.
The new insertion points are $\{\}$.

This iteration adds no future ideals.

From here onward we no longer list the future ideals as they are
generated as if there are any, there are too many of them to list
conveniently.

\paragraph{The ideal \Ideal{A,C}} \disconnect

\Figure{\ipeFigACplusBD}

The deletion points are $\{C\}$.
The number of insertion points is 3.
The existing insertion points are $\{B,D\}$.
The new insertion points are $\{G\}$.

\Figure{\ipeFigACplusBDG}

\paragraph{The ideal \Ideal{A,D}} \disconnect

\Figure{\ipeFigADplusBC}

The deletion points are $\{D\}$.
The number of insertion points is 3.
The existing insertion points are $\{B,C\}$.
The new insertion points are $\{H\}$.

\Figure{\ipeFigADplusBCH}

\paragraph{The ideal \Ideal{A,C,D}} \disconnect

\Figure{\ipeFigACDplusBGH}

The deletion points are $\{C,D\}$.
The number of insertion points is 4.
The existing insertion points are $\{B,G,H\}$.
The new insertion points are $\{I\}$.

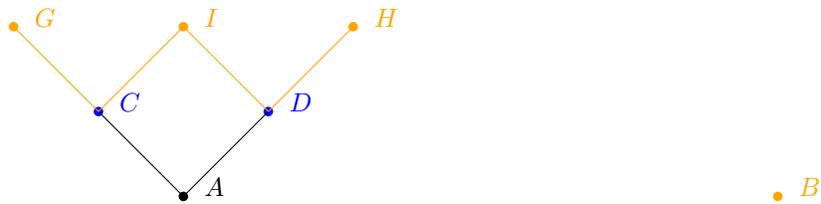
\Figure{\ipeFigACDplusBGHI}

\paragraph{The ideal \Ideal{A,B,C}} \disconnect

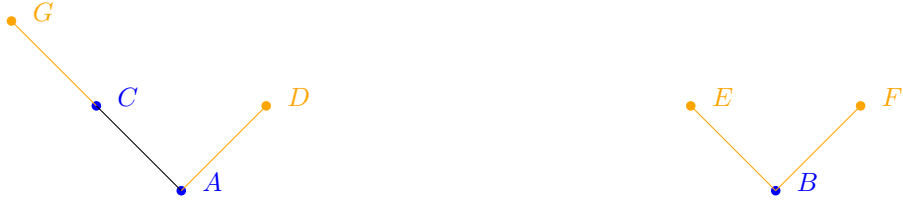
\Figure{\ipeFigABCplusDEFG}

The deletion points are $\{B,C\}$.
The number of insertion points is 4.
The existing insertion points are $\{D,E,F,G\}$.
The new insertion points are $\{\}$.

\vspace{\baselineskip}

\emph{Ideals that split into two non-empty disjoint parts (i.e., that
contain both $A$ and $B$) never have new insertion points.}

\paragraph{The ideal \Ideal{A,B,D}} \disconnect

No new insertion points.

\paragraph{The ideal \Ideal{A,B,C,D}} \disconnect

No new insertion points.

\paragraph{The ideal \Ideal{B,E}} \disconnect

\Figure{\ipeFigBEplusAF}

The deletion points are $\{E\}$.
The number of insertion points is 3.
The existing insertion points are $\{A,F\}$.
The new insertion points are $\{J\}$.

\Figure{\ipeFigBEplusAFJ}

\paragraph{The ideal \Ideal{B,F}} \disconnect

\Figure{\ipeFigBFplusAE}

The deletion points are $\{F\}$.
The number of insertion points is 3.
The existing insertion points are $\{A,E\}$.
The new insertion points are $\{K\}$.

\Figure{\ipeFigBFplusAEK}

\paragraph{The ideal \Ideal{B,E,F}} \disconnect

\Figure{\ipeFigBEFplusAJK}

The deletion points are $\{E,F\}$.
The number of insertion points is 4.
The existing insertion points are $\{J,K\}$.
The new insertion points are $\{L\}$.

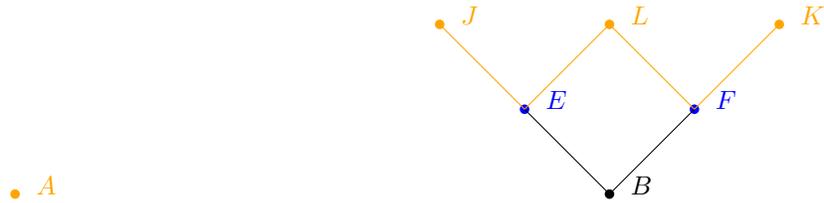
\Figure{\ipeFigBEFplusAJKL}

The process continues to construct the poset of points, which is two
disjoint copies of the quadrant.

This example generalizes to show that there is exactly one
$d$-differential distributive lattice for every positive integer $d$,
that is, with each point having
the weight 1 and each lattice element having the differential degree $d$.
This lattice is $\mbbY^{\times d}$.

\section{Proving there are no triple orphans}

The process can be analyzed to prove other facts.  For example, we can
prove:  \emph{If the differential degree $d$ is constant over the
lattice, no point is covered by three orphans,} where we define
\emph{orphan} as a point that covers exactly one point.

Consider a point $P$ which is assumed to be covered by at least three
orphans.  Consider the principal ideal $[P]$ generated by
$P$, the set of all points $\leq P$:

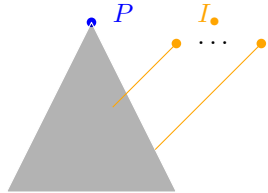
\Figure{\ipeFigP}

By construction, the only deletion point of this ideal is $P$.
Since this is the first ideal containing $P$, no old insertion points
cover $P$, but there is a (possibly empty) set of insertion points
$I_\bullet$ that cover elements less than $P$.
We label the new insertion points $A_\bullet$, using $\bullet$ as the
free index in this notation.  The $A_\bullet$ necessarily cover
$P$ and only $P$, that is, they are the orphans covering $P$.
The differential condition requires
\begin{equation}
  P + d = \sum I_\bullet + \sum A_\bullet.  \label{eq:1}
\end{equation}
(We use the name of a point to also denote its weight.
In summations, we use $\bullet$ as the bound index.)
Each $A_\bullet$ covers exactly $P$.

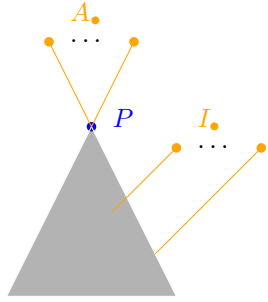
\Figure{\ipeFigPA}

Consider the principal ideal $[A_i]$:

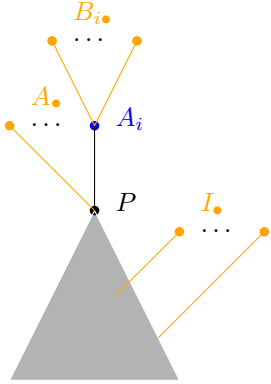
\Figure{\ipeFigPAiB}

The only deletion point is $A_i$.
The old insertion points are the $I_\bullet$ and the $A_\bullet$
except $A_i$.
We label the new insertion points
$B_{i \bullet}$.  They necessarily cover $A_i$ and only $A_i$.
The differential condition requires
\begin{align}
A_i + d & = \sum I_\bullet + \sum A_\bullet - A_i + \sum B_{i \bullet} \\
  2A_i + d & = \sum I_\bullet + \sum A_\bullet + \sum B_{i \bullet} \\
\interject{Substituting (\ref{eq:1}),}
2A_i - P & = \sum B_{i \bullet}.  \label{eq:4}
\end{align}

Consider the ideal $[A_i,A_j]$ which is generated by $A_i$ and $A_j$:

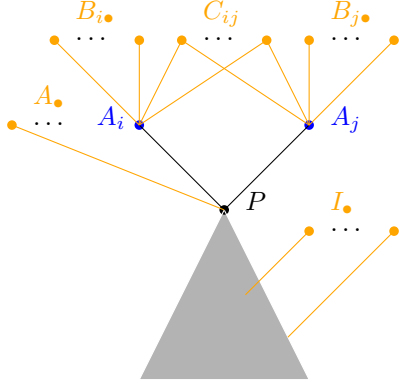
\Figure{\ipeFigPAiAjBC}

The only deletion points are $A_i$ and $A_j$.
The old insertion points are the $I_\bullet$, the $A_\bullet$
except $A_i$ and $A_j$, the $B_{i \bullet}$, and the $B_{j \bullet}$.
We label the new insertion points
$C_{ij \bullet}$.  They necessarily cover $A_i$ and $A_j$ and no other points.
The differential condition requires
\begin{align}
A_i + A_j + d & =
  \sum I_\bullet + \sum A_\bullet - A_i - A_j
  + \sum B_{i \bullet} + \sum B_{j \bullet} + \sum C_{ij \bullet} \label{eq:5} \\
2A_i + 2A_j + d & =
  \sum I_\bullet + \sum A_\bullet
  + \sum B_{i \bullet} + \sum B_{j \bullet} + \sum C_{ij \bullet} \label{eq:6} \\
\interject{Substituting (\ref{eq:4}) twice,}
2A_i + 2A_j + d & =
  \sum I_\bullet + \sum A_\bullet
  + 2A_i - P + 2A_j -P + \sum C_{ij \bullet}, \label{eq:7} \\
\interject{Substituting (\ref{eq:1}),}
2A_i + 2A_j + d & =
  P + d
  + 2A_i - P + 2A_j -P + \sum C_{ij \bullet} \\
P & = \sum C_{ij \bullet}. \label{eq:11}
\end{align}

Consider the ideal $[A_i,A_j,A_k]$ which is generated by $A_i$, $A_j$,
and $A_k$.  This ideal is too messy to draw easily, but we can analyze
it the same way as the preceding two ideals:

The only deletion points are $A_i$, $A_j$ and $A_k$.
The old insertion points are the $I_\bullet$; the $A_\bullet$
except $A_i$, $A_j$, and $A_k$; the $B_{i \bullet}$, the
$B_{j \bullet}$, and the $B_{k \bullet}$; and
the $C_{ij \bullet}$, the $C_{ik \bullet}$, and the $C_{jk \bullet}$.
We label the new insertion points
$D_{ijk \bullet}$.  They necessarily cover $A_i$, $A_j$, and $A_k$
and no other points.
The differential condition requires
\begin{align}
  A_i + A_j + A_k + d & =
  \sum I_\bullet + \sum A_\bullet - A_i - A_j - A_k
  + \sum B_{i \bullet} + \sum B_{j \bullet} + \sum B_{k \bullet} \notag \\
  & \hspace{3em} + \sum C_{ij \bullet} +  \sum C_{ik \bullet} +   \sum C_{jk \bullet} +
  \sum D_{ijk \bullet} \\
2A_i + 2A_j + 2A_k + d & =
  \sum I_\bullet + \sum A_\bullet
  + \sum B_{i \bullet} + \sum B_{j \bullet} + \sum B_{k \bullet} \notag \\
  & \hspace{3em} + \sum C_{ij \bullet} +  \sum C_{ik \bullet} +   \sum C_{jk \bullet}
  + \sum D_{ijk \bullet} \\
\interject{Substituting (\ref{eq:4}) and (\ref{eq:11}) three times,}
2A_i + 2A_j + 2A_k + d & =
  \sum I_\bullet + \sum A_\bullet
  + 2A_i - P + 2A_j - P + 2A_k - P \notag \\
  & \hspace{3em} + P + P + P
  + \sum D_{ijk \bullet} \\
\interject{Substituting (\ref{eq:1}),}
2A_i + 2A_j + 2A_k + d & =
  P + d
  + 2A_i - P + 2A_j - P + 2A_k - P
  + P + P + P
  + \sum D_{ijk \bullet} \\
0 & = P + \sum D_{ijk \bullet} \label{eq:20}
\end{align}
But since $P$ and all of the $D_{ijk \bullet}$ are positive,
(\ref{eq:20}) is impossible, so there is no $P$ covered by three orphans.

\vspace{3em}

\end{document}